\documentclass[11pt,a4paper]{article}
\usepackage{amsmath, ifthen}
\usepackage{amssymb, amsfonts}

\title{Effect of linear lumping on \\ controllability and observability}
\author {Zs\'{o}fia Horv\'{a}th\thanks{Department of Mathematical Analysis,
Budapest University of Technology and Economics,
Budapest, H-1111 Egry J. u. 1., HUNGARY
}}
\newcommand{\Ref}[1]{(\ref{#1})}
\newcommand{\R}{\mathbb{R}}
\newcommand{\N}{\mathbb{N}}
\newcommand{\rank}{\mathrm{rank}}

\newcommand{\lS}{linear system}
\newcommand{\laS}{lumped system}
\newcommand{\dE}{differential equation}
\newcommand{\cC}{completely controllable}
\newcommand{\cO}{completely observable}

\renewcommand{\emph}{\textbf}

\newcommand{\Span}{\mathrm{Span}}

\newcommand{\torol}[1]{}

\newtheorem{elnev}{Definition}
\newtheorem{tetel}{Theorem}
\newtheorem{megj}{Remark}

\begin{document}
\maketitle
\vfill{\Large
\begin{description}
\item[Email address:] hzsofi@math.bme.hu
\item[Phone:] 361-463-2324
\item[Fax:] 361-463-3172
\item[Running head:]
Lumping, controllability and observability
\end{description}
}
\eject

\begin{abstract}
The effect of linear lumping, linear transformation to reduce the
number of state variables on controllability and observability of
linear differential equations has been studied. Controllability of
the original system implies the controllability of the lumped
system. Examples taken from reaction kinetics illustrate our
results.
\subsection*{KEY WORDS}
completely controllable, lumping matrix, compartmental system,
M-matrix
\subsection*{MSC}
80A30, 93B17, 34A30
\end{abstract}

\section{Introduction}

Dealing with modeling of real questions the large number of
variables is generally a problem. To have a model which can easily
be treated, one possible way is to reduce the number of variables
by a method called lumping \cite{lump}. By controllability of a
system we mean that it can be brought from any position to any
other position in a finite amount of time. On the other hand,
observability of a system means that we can determine the initial
state of the system from the knowledge of an input-output pair
over a certain period of time.

In this paper we study the effect
of linear lumping on such properties of the system as
controllability and observability and apply the results to
compartmental systems. Our main result is that a controllable
compartmental system remains controllable under lumping.

\section{Basic notions of lumping and control theory}

In this section we collect the basis of the mathematical theory of
controllability \cite[page 16]{Chen}, observability \cite[page
26]{Chen} and lumping \cite[page 1534]{lump}. Before turning to
the formal definitions we mention that controllability means that
any prescribed concentration can be attained using an appropriate
control input, observability has the meaning that one can
reconstruct the history of the process when knowing the present
concentration composition.

Let $n,r,p\in{\N}, A\in{\R^{n\times n}}, B\in{\R^{n\times r}},
C\in{\R^{p\times n}}$ and let us investigate:
\begin{eqnarray}
\dot x(t)&=&Ax(t)+Bu(t)
\label{rsz}
\\  y(t)&=&Cx(t)
\label{mf}
\end{eqnarray}
with $x(t)$ in $R^n$ as the dependent variable of the linear
differential equations at time $t$ in $R$ and $u(t)$ in $R^r$
denoting the  bounded and measurable control function and $ y $
with values in $R^p$ the observation function.

The fundamental definitions and the Kalman rank conditions below
can be found in any textbook on linear system theory, e.g. in
\cite{Chen}.

\begin{elnev}
A \lS\ with a state-space description given by \Ref{rsz} is said
to be \emph{\cC} if, starting from any position $x_{0}$ in $R^n$
at any initial time $t_{0}$, the state vector $x$ can be brought
to any other position $x_{1}$ in $\R^n$ in a finite amount of time
by some control function $u$.
\end{elnev}

\begin{tetel} \label{ir}
A \lS\ described by \Ref{rsz} is \cC\ if and only if the $n\times rn$
 matrix $W_{AB}:=[B \mid AB \mid \ldots \mid A^{n-1}B]$ has rank $n$.
\end{tetel}

Let $t_{0},t_{1}\in{\R}; t_{0}<t_{1}$.

\begin{elnev}
A \lS\ with the state-space description \Ref{rsz}--\Ref{mf} has
the observability property on an interval $(t_{0},t_{1})$, if any
input-output pair $(u(t),y(t))$, $t_{0}\leq t\leq t_{1}$, uniquely
determines the initial state $x(t_{0})$. Furthermore
\Ref{rsz}--\Ref{mf} is said to be observable at an initial time
$t_{0}$ if it has the observability property on some interval
$(t_{0},t_{1})$ where $t_{1}>t_{0}$. It is said to be \emph{\cO}
if it is observable at every initial time $t_{0}$.
\end{elnev}

\begin{tetel} \label{megf}
A \lS\ described by \Ref{rsz}--\Ref{mf}  is \cO\ if and only if
the
 $n\times pn$ matrix $V_{CA}:=[C^{\top} \mid A^{\top}C^{\top} \mid \ldots \mid (A^{\top})^{n-1}
C^{\top}]$ has rank $n$.
\end{tetel}

Let us consider
\begin{eqnarray}
\dot x(t)&=&A^{\top}x(t)+C^{\top}v(t)
\label{drsz}
\\  z(t)&=&B^{\top}x(t)
\label{dmf}
\end{eqnarray}
with $v(t)$ in $R^p$ denoting the  bounded and measurable control
function and $z(t)$ in $R^r$ the observation function.

\begin{tetel} \label{du}
The \lS s \Ref{rsz}--\Ref{mf} and \Ref{drsz}--\Ref{dmf}
respectively described above are \emph{dual to each other} in the
sense that \Ref{rsz}--\Ref{mf} is \cO\ if and only if \Ref{drsz}
is \cC\, and \Ref{rsz} is \cC\ if and only if \Ref{drsz},
\Ref{dmf} is \cO\ {\em \cite[page 32]{Chen}}.
\end{tetel}

Let us introduce a slightly modified notion of lumpability
\cite[page 1262]{Mmx}.

\begin{elnev}
Suppose $l \in{\N},$ $l \leq n,$ and $M \in{\R^{l \times n}},$
$\rank(M)=l$. If for all solutions $x(t)$ to \Ref{rsz}
\begin{eqnarray}
\hat x(t):=Mx(t)
\label{re}
\end{eqnarray}
obeys a differential equation
\begin{eqnarray}
\dot{\hat x}(t)=\hat A \hat x(t)+ \hat u(t) \label{rersz}
\end{eqnarray}
with $\hat A \in{\R^{l \times l}}$, and \quad $\hat u(t):=MBu(t)$
then \Ref{rsz} is said to be \emph{exactly lumpable} to
\Ref{rersz} by $M$. The pair of matrices $M$ and $\hat A$ is
sometimes referred to as a lumping scheme.
\end{elnev}

\begin{megj} \label{Akalap}
{\em It can be shown that \Ref{rsz} is exactly lumpable to
\Ref{rersz} by $M$ if and only if there exists $\hat A \in{\R^{l
\times l}},$ such that ${\hat A}M=MA$, so we get the lumping
scheme consisting of $M$ and $\hat A=MAM^{\top}(MM^{\top})^{-1}$
\cite{lump}. We note that the inverse $(MM^{\top})^{-1} \in{\R^{l
\times l}}$ exists, since $\rank(M)=l$ and
$M^{\top}(MM^{\top})^{-1}$ is a generalized inverse of $M$
satisfying $MM^{\top}(MM^{\top})^{-1}=I_{l \times l}$ with $I_{l
\times l}$ being the $l$-identity matrix. It can be proved that
there exists a matrix $\hat A$ such that ${\hat A}M=MA$ if
$$M=N\left[ \begin{array}{c} f_{1}^{\top} \\ \vdots \\ f_{l}^{\top}\end{array} \right]$$
where $f_{i}$ $(i=1, \ldots , l)$ are any independent, real
eigenvectors of the matrix $A$ and the matrix $N \in{\R^{l\times
l}}$ is nonsingular \cite[page 117]{inv}. }
\end{megj}

\torol{
\begin{megj} \label{Jinv}
{\em If the geometric multiplicity of real eigenvalues of $A$
equals to 1 and nonreal eigenvalues appear in complex conjugate
pairs and there exists a real similarity transformation matrix $F
\in{\R^{n \times n}}$ such that $A^{\top}F=FJ$ with $J
\in{\R^{n\times n}}$ being the Jordan form of $A$ then the lumping
matrix can be obtained as
$$M=N(FI)^{\top}$$ where $I \in{\R^{n
\times l}}$ and $ I= \{e_{1},\ldots, e_{l} \},$ its columns
$e_{i}$ \quad $i=1,\ldots,l$ are the unit vector of $\R^n$ with 1
as its $i$th element, and 0 for the rest of the elements, and
matrix $N \in{\R^{l\times l}}$ is nonsingular. Using the fact that
$\hat A=MAM^{\top}(MM^{\top})^{-1}$ by Remark \ref{Akalap} and
note that there exists the matrix $(MM^{\top})^{-1}$ since
$\rank(M)=l,$ we get $\hat
A=N(FI)^{\top}A(FI)((FI)^{\top}(FI))^{-1}N^{-1}$.}
\end{megj}
}%torolvege

\section{Fundamentals}
A typical chemical reaction step is written schematically as
\begin{equation}
 Reactant {\rightarrow} Product.
\end{equation}
The reactant and product are each called a complex, which is a set
of elements with associated coefficients. The elements that make
up complexes are called species, and can be anything that
participates in a reaction, typically a chemical element, a
molecule, or an electron. The species that are on the left side of
the equation are used up, and those on the right are created when
the reaction step occurs. The coefficient that a species takes
indicates what proportion of it is created or used in the
reaction, and by convention it is always a non-negative integer.
The goal of chemical reaction theory is to monitor how the
concentration of each species changes over time \cite[page
1539]{lump}. The fundamental constituent of a chemical reaction is
the species, the concentration of which we are interested in
monitoring. A complex is a sum of species with integer
coefficients, and a reaction is a pair of complexes with an
ordering (to distinguish between products and reactants). Also a
subclass of the class of chemical processes can be modeled as
compartmental systems \cite{Brochot}, \cite{compart}. A
compartmental system consists of several compartments with more or
less homogeneous amounts of material. The compartments interact by
processes of transport and diffusion. The dynamics of a
compartmental system is derived from mass balance considerations.
The mathematical theory of compartmental systems is of major
importance to pharmacokineticists, physiologists \cite{compart}.
Sometimes it is useful to reduce a model to get a new one with a
lower dimension. The technique's name is lumping, i.e. reduction
of the number of variables by grouping them via a linear or
nonlinear function. Very often problems in the physical, and
social sciences can be reduced to problems involving matrices
which, due to certain constraints, have some special structure
\cite[page 132]{inv}.

%torolvege

\section{Effect of lumping on controllability and observability}

\torol{
 Let $A \doteq A_{n}=-s(I_{n \times n}-T_{n \times n}),
\quad s>0.$
\begin{elnev}
$A_{n} \in{\R^{n\times n}}$ can then be expressed in the form
\begin{equation} \label{felbontas}
A_{n}=-s(I_{n \times n}-T_{n \times n}), \quad s>0,
\end{equation}
$s\in{\R},$ $T_{n \times n}\in{\R^{n\times n}}$ nonnegative,
irreducible and symmetric matrix. Any matrix $A_{n}$ of the form
\Ref{felbontas} for which $ \varrho(T_{n\times n})\leq 1$ (the
spectral radius of $\quad$ $T_{n \times n}:$ $\varrho(T_{n \times
n})= \max_{h=1,2, \ldots , n} \mid \lambda_{h} \mid ,$
$\lambda_{h} \in{\R}$ there is eigenvalues of $T_{n \times n}$) is
called an $M-matrix.$
\end{elnev}
}%torolvege

Our main result is expressed in the following statement.

\begin{tetel} \label{redir}
Let us assume that \Ref{rsz} is exactly lumpable to \Ref{rersz} by
$M$ and
 the \lS\ \Ref{rsz} is \cC, then the \lS\ \Ref{rersz} is
also \cC .
\end{tetel}
\textit {Proof.} Let \Ref{rsz} be \cC\, then according to Theorem \ref{ir} we know
that $\rank(W_{AB})=n$. Furthermore, using the fact that
$$\hat AM=MA,\ldots,\hat A^{l-1}M=MA^{l-1}$$
we get
$$\hat W_{AB}:=[MB \mid {\hat A}MB \mid \ldots \mid (\hat A)^{l-1}MB]=M [B \mid AB \mid
\ldots \mid A^{l-1}B].$$ This implies that \Ref{rersz} is \cC\ if
and only if the $l\times rl$ matrix
$$M [B \mid AB \mid \ldots \mid A^{l-1}B]$$
has rank $l$ by Theorem \ref{ir}. Let us assume that the rank of
$\hat W_{AB}$ is less than $l$, then there is a nonzero vector $b
\in{\R^l}$ with $b^{\top}\hat W_{AB}=0 \in{\R^{rl}}$ or
equivalently $b^{\top}MB=b^{\top}\hat AMB= \ldots =b^{\top}(\hat
A)^{l-1}MB=0 \in{\R^r}$. Application of the Cayley--Hamilton
Theorem now gives
$$(\hat A)^{l}=\gamma_{1}(\hat A)^{l-1}+\gamma_{2}(\hat A)^{l-2}+\ldots+\gamma_{l}I_{l \times l}$$
where $I_{l \times l}$ represents the $l\times l$ unit matrix and
$\gamma_{1},\ldots,\gamma_{l}$ are suitable constants, so
$b^{\top}(\hat A)^{l}MB=0 \in{\R^r}$. With induction we can derive
that $b^{\top}(\hat A)^{l+j}MB=0 \in{\R^r}$ for
$(j=1,2,\ldots,n-1-l)$ also, thus
$$b^{\top}[MB \mid {\hat A}MB \mid \ldots \mid (\hat A)^{n-1}MB]=b^{\top}M [B \mid AB \mid
\ldots \mid A^{n-1}B]=0 \in{\R^{rn}}.$$ Since $\rank(M)=l$,
therefore $b^{\top}M \neq0 \in{\R^n}$ thus we get that the rows of
matrix $W_{AB}=[B \mid AB \mid \ldots \mid A^{n-1}B]$ are linearly
dependent, which is a contradiction. Consequently matrix $\hat
W_{AB}$ has to have full rank and hence \Ref{rersz} is \cC .

\begin{megj} \label{redirm}
{\em The complete controllability of system \Ref{rersz} does not
imply the complete controllability of system \Ref{rsz}. Since even
if the rank of $\hat W_{AB}$ is $l$, the rank of $W_{AB}$ can be
less then $n$. A concrete example will also be given below on page
\pageref{pelda1}.}
\end{megj}

\begin{megj} \label{egyen}

{\em In the case of $l=n$ we get that $M$ is a non-singular
$n\times n$ matrix, thus
\begin{eqnarray}
\rank(\hat W_{AB})=&\rank(M [B \mid AB \mid \ldots \mid
A^{l-1}B])&= \nonumber
\\=&\rank([B \mid AB \mid \ldots \mid A^{n-1}B])&=\rank(W_{AB}).
\nonumber
\end{eqnarray}
Therefore if \Ref{rsz} is \cC\ then \Ref{rersz} is also \cC\
 and vice versa. }
\end{megj}

\begin{megj} \label{drem}
{\em Let us assume that \Ref{rsz} is exactly lumpable to
\Ref{rersz} by $M$ then \Ref{drsz} is exactly lumpable to
\begin{eqnarray}
\dot{\hat x}(t)=\tilde A \hat x(t)+\ MC^{\top}v(t)
\label{drersz}
\end{eqnarray}
by $ M,$ where $\tilde A \in{\R^{l \times l}}$ such that $
MA^{\top}=\tilde AM,$ so $\tilde
A=MA^{\top}M^{\top}(MM^{\top})^{-1}. $ If we assume that
\Ref{rsz}, \Ref{mf} is \cO\ then \Ref{drsz} is \cC\ by Theorem
\ref{du}, so \Ref{drersz} is also \cC\ by the application of
Theorem \ref{redir}.}
\end{megj}

\begin{megj} \label{remf}
{\em The complete controllability of system \Ref{drersz} does not
imply the complete observability of system \Ref{rsz}--\Ref{mf}.
The controllability matrix of system \Ref{drersz} is $[MC^{\top}
\mid \tilde AMC^{\top} \mid \ldots \mid (\tilde A)^{l-1}
MC^{\top}]=M[C^{\top} \mid A^{\top}C^{\top} \mid \ldots \mid
(A^{\top})^{l-1} C^{\top}]$ where we have used $\tilde
AM=MA^{\top},\ldots,\tilde A^{l-1}M=M(A^{\top})^{l-1}.$ Even if it
is of has rank $l$, the rank of $V_{CA}$ can be less then $n$. A
concrete example will also be given below on page
\pageref{pelda2}.}
\end{megj}

%torolvege
\section{Examples}

In this section we will give examples to illustrate our results.

Let us consider the following chemical reaction, a special case of
a compartmental system \cite[page 3]{E.T.}:
\begin{equation} \label{reakcio}
{\cal X}_{1} \stackrel{k}{\rightarrow}  {\cal X}_{2} \quad {\cal
X}_{2} \stackrel{k}{\rightarrow} {\cal X}_{1} \quad {\cal X}_{2}
\stackrel{k}{\rightarrow} {\cal X}_{3} \quad {\cal X}_{3}
\stackrel{k}{\rightarrow} {\cal X}_{2}
\end{equation}
where ${\cal X}_{i}$ \quad $(i=1,2,3)$ is the $i$th chemical
component or species and the positive real number $k$ is the
uniform reaction rate constant. The example may be degenerate from
the point of view of kinetics, however it may be considered as a
three stage approximation of diffusion in a tube. We are
interested in the time evolution of the quantities of chemical
components. If we assume that the physical circumstances are
ideal, i.e. in the given reaction the temperature, the pressure,
and the volume of the vessel are constant we can build up the mass
action-type model of the reaction \Ref{reakcio}:
\begin{eqnarray}
\dot x_{1}(t)&=&-kx_{1}(t)+kx_{2}(t)
\nonumber
\\ \dot x_{2}(t)&=&kx_{1}(t)-2kx_{2}(t)+kx_{3}(t)
\label{reakcegy}
\\ \dot x_{3}(t)&=&kx_{2}(t)-kx_{3}(t)
\nonumber
\end{eqnarray}
where $x_{i}(t)$ \quad $(i=1,2,3)$ is to be interpreted as the
concentration of the species \( {\cal X}_{i} \) at time $t$.
Equation \Ref{reakcegy} is said to be the induced kinetic \dE\ of
\Ref{reakcio}.

To construct the lumping matrix $M$ we use the following facts: if
every element of the matrix $M=N\left[
\begin{array}{c} f_{1}^{\top} \\ \vdots \\ f_{l}^{\top}\end{array}
\right]$ is nonnegative, moreover for every row of $M$ there
exists an element from that row which is the only nonzero element
of its column, then the lumped system of the induced kinetic
differential equation of a reaction is also a induced kinetic \dE\
of a reaction \cite{F.Gy.}.

Using this and the fact that the independent, real eigenvectors of
$$A=\left[ \begin{array}{rrr} -k&k&0 \\ k&-2k&k \\ 0&k&-k \end{array} \right]$$
are
$$f_{1}^{\top}=[1,1,1], \quad f_{2}^{\top}=[1,0,-1], \quad f_{3}^{\top}=[1,-2,1] $$
set
$$M:=\left[ \begin{array}{rr} 1&1 \\ 1&-1 \end{array} \right]
\left[ \begin{array}{rrr} 1&1&1\\ 1&0&-1 \end{array} \right]=
\left[ \begin{array}{rrr} 2&1&0\\ 0&1&2 \end{array} \right]$$ in
all the examples of this section. In this case the new variables
are
$$\left[ \begin{array}{c} \hat x_{1} \\ \hat x_{2} \end{array} \right]=
\left[ \begin{array}{rrr} 2&1&0\\ 0&1&2 \end{array} \right] \left[
\begin{array}{c}  x_{1} \\ x_{2} \\ x_{3} \end{array} \right]=
\left[ \begin{array}{c}  2x_{1}+x_{2} \\ x_{2}+2x_{3} \end{array}
\right]$$ and the lumped system is
\begin{eqnarray}
\dot{\hat x}_{1}&=& - \frac{k}{2} \hat x_{1}+ \frac{k}{2} \hat x_{2}
\nonumber
\\ \dot{\hat x}_{2}&=& \frac{k}{2} \hat x_{1}- \frac{k}{2} \hat x_{2}
\nonumber
\end{eqnarray}
which is the induced kinetic \dE\ of the reaction

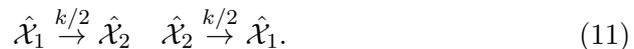
\begin{equation}
{\cal{\hat X}}_{1} \stackrel{k/2}{\rightarrow} {\cal{\hat X}}_{2}
\quad {\cal{\hat X}}_{2} \stackrel{k/2}{\rightarrow} {\cal{\hat
X}}_{1}. \label{reakcior}
\end{equation}

Let us remark that the new variables can be considered as groups
of the old ones measured together, since they are nonnegative
linear combinations of the old ones.

We shall use the notion and properties of an M-matrix below.

\begin{elnev}
Let $n\in{\N}.$ If the matrix $A_{n} \in{\R^{n\times n}}$ can be
expressed in the form
\begin{equation} \label{felbontas}
A_{n}=-s(I_{n \times n}-T_{n \times n}),
\end{equation}
 where $s\in{\R^+},$ $T_{n \times n}\in{\R^{n\times n}}$ is a nonnegative,
irreducible and symmetric matrix with the spectral radius
$\varrho$ obeying
$$
\varrho(T_{n \times n})= \max_{h=1,2, \ldots ,
n} \mid \lambda_{h} \mid \leq 1,
$$
where $\lambda_{h} \in{\R}$ are
eigenvalues of $T_{n \times n},$ then $A_{n}$ is an M-matrix.
\end{elnev}

Since both the coefficient matrix of the original system and that
of the lumped one are of the form \Ref{felbontas}:
$$A=\left[ \begin{array}{rrr} -k&k&0 \\ k&-2k&k \\
0&k&-k
\end{array} \right]=-2k \left(\left[ \begin{array}{rrr} 1&0&0 \\ 0&1&0 \\ 0&0&1 \end{array}
\right]- \left[ \begin{array}{rrr} \frac{1}{2}&\frac{1}{2}&0 \\ \frac{1}{2}&0&\frac{1}{2} \\
0&\frac{1}{2}&\frac{1}{2}
\end{array} \right] \right),$$
$$ \hat A =\left[ \begin{array}{rr} -k&k \\
k&-k \end{array} \right]=-2k \left(\left[ \begin{array}{rr} 1&0 \\
0&1 \end{array} \right]-\left[ \begin{array}{rr} \frac{1}{2}&\frac{1}{2} \\
\frac{1}{2}&\frac{1}{2} \end{array} \right]  \right),$$ where
$k>0$ and $\varrho(T_{3 \times 3})= \max [ \mid \frac{1}{2} \mid ,
\mid 1 \mid , \mid -\frac{1}{2} \mid ]=1,$ and $\varrho(T_{2
\times 2})= \max [ \mid 0 \mid , \mid 1 \mid]=1,$ thus both
matrices is an M-matrix. A similar construction can be obtained
for an even number of species as follows.

%torolvege

First, let $n=4 \theta +2 \quad (\theta=0,1, \ldots ) $ and let us
consider the following chemical reaction:
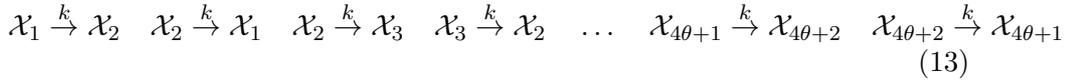
\begin{equation} \label{reakcio8} {\cal X}_{1}
\stackrel{k}{\rightarrow}  {\cal X}_{2} \quad {\cal X}_{2}
\stackrel{k}{\rightarrow} {\cal X}_{1} \quad {\cal X}_{2}
\stackrel{k}{\rightarrow} {\cal X}_{3} \quad {\cal X}_{3}
\stackrel{k}{\rightarrow} {\cal X}_{2}\quad \ldots \quad {\cal
X}_{4 \theta +1} \stackrel{k}{\rightarrow} {\cal X}_{4 \theta+2}
\quad {\cal X}_{4 \theta+2} \stackrel{k}{\rightarrow} {\cal X}_{4
\theta +1}
\end{equation}
then two independent, real eigenvectors of the matrix
$$A_{4 \theta+2}= \left[ \begin{array}{rrrrrr} -k&k&0& \ldots &0 &0 \\ k&-2k&k& \ldots &0 &0 \\
 0&k&-2k& \ldots &0 &0 \\ \vdots& \vdots & \vdots & \ddots & \vdots & \vdots
 \\ 0&0&0& \ldots &-2k &k \\ 0&0&0& \ldots &k &-k \end{array} \right]$$
are
$$f_{1}^{\top}=[1,1,1,1, \ldots ,1,1], \quad f_{2}^{\top}=[-1,1,1,-1, \ldots ,-1,1] $$
so if we take
$$M:=\left[ \begin{array}{rr} 1&1 \\ 1&-1 \end{array} \right]
\left[ \begin{array}{rrrrrrr} 1&1&1&1& \ldots &1&1\\
-1&1&1&-1 & \ldots &-1&1
\end{array} \right]= \left[ \begin{array}{rrrrrrr} 0&2&2&0& \ldots &0&2\\
2&0&0&2& \ldots &2&0
\end{array} \right],$$
then we obtain
 $$ A_{2}=\left[ \begin{array}{rr} -k&k \\ k&-k \end{array}
\right]$$ and the lumped system is again \Ref{reakcior}.

Second, let $n=4 \theta \quad (\theta=1,2, \ldots ) ,$ then two
independent, real eigenvectors of the matrix $A_{4 \theta}$ are
$$f_{1}^{\top}=[1,1,1,1, \ldots ,1,1 ], \quad f_{2}^{\top}=[-1,1,1,-1, \ldots ,1,-1] $$
so if we take
$$M:=\left[ \begin{array}{rr} 1&1 \\ 1&-1 \end{array} \right]
\left[ \begin{array}{rrrrrrr} 1&1&1&1& \ldots &1&1 \\
-1&1&1&-1 & \ldots &1&-1
\end{array} \right]= \left[ \begin{array}{rrrrrrr} 0&2&2&0& \ldots & 2&0\\
2&0&0&2& \ldots &0&2
\end{array} \right],$$
then we get $ A_{2}$ and the lumped system is \Ref{reakcior}. We
note that $A_{n},$ $n=4 \theta ,\quad or \quad n=4 \theta +2 \quad
(\theta=1,2, \ldots),$ is again of the form \Ref{felbontas},
with $ s=2k>0$ and
$$T_{n\times n}= \left[ \begin{array}{rrrrrr} \frac{1}{2}&\frac{1}{2}&0& \ldots &0 &0 \\
 \frac{1}{2}&0&\frac{1}{2}& \ldots &0 &0 \\
 0&\frac{1}{2}&0& \ldots &0 &0 \\ \vdots& \vdots & \vdots & \ddots & \vdots & \vdots
 \\ 0&0&0& \ldots &0 &\frac{1}{2} \\ 0&0&0& \ldots &\frac{1}{2} &\frac{1}{2} \end{array}
 \right]>0,$$ and $\lambda_{h}= \cos(\frac{(h-1) \pi}{n}),$ \quad $h=1,2 \ldots , n$
 are eigenvalues of $T_{n \times n}$ \cite[page 177]{inv}: and
$\varrho(T_{n \times n})= \mid \cos(\frac{0 \pi}{n}) \mid =1$ is
the spectral radius of $T_{n \times n},$ thus $A_{n}$ in an
M-matrix in this case too. Furthermore, the eigenvalues of
$A_{n}$ are given by $-2k(1-\lambda_{h})=-4k (\sin(\frac{(h-1)
\pi}{2n}))^{2},$ \quad $h=1,2 \ldots , n,$ so they are real and
negative, and zero, and any independent, real eigenvectors of
$A_{n}$ are $f_{1q}=\frac{1}{\sqrt{n}}$; $f_{\eta
q}=\sqrt{\frac{2}{n}} \cos(\frac{(2q-1)(\eta-1) \pi}{2n})$ \quad
$\eta=2,3 \ldots , n,$ \quad $q=1,2 \ldots , n$.

\subsection{Examples of the effect of lumping on controllability}
\subsubsection{Controllable original---controllable lumped system}
First, let us consider the \dE
$$\frac{d}{dt} \left[ \begin{array}{c}  x_{1} \\ x_{2} \\ x_{3} \end{array} \right]=
\left[ \begin{array}{rrr} -k&k&0 \\ k&-2k&k \\ 0&k&-k \end{array} \right]
\left[ \begin{array}{c}  x_{1} \\ x_{2} \\ x_{3} \end{array} \right]+
\left[ \begin{array}{rrr} 1&0&0 \\ 0&1&0 \\ 0&0&1 \end{array} \right]
\left[ \begin{array}{c}  u_{1} \\ u_{2} \\ u_{3} \end{array} \right].$$
then
$$W_{AB}=[B \mid AB \mid A^{2}B]=
\\ \left[ \begin{array}{rrrrrrrrr} 1 &0 &0 &-k &k &0 &2k^2 &-3k^2 &k^2
\\ 0 &1 &0 &k &-2k &k &-3k^2 &6k^2 &-3k^2 \\ 0 &0 &1 &0 &k &-k &k^2 &-3k^2 &2k^2
\end{array} \right].$$
This matrix has rank 3, therefore this system is \cC .
The \laS\ is
$$\frac{d}{dt} \left[ \begin{array}{c} \hat x_{1} \\ \hat x_{2} \end{array} \right]=
\left[ \begin{array}{rr} -\frac{k}{2}&\frac{k}{2} \\
\frac{k}{2}&-\frac{k}{2} \end{array} \right] \left[
\begin{array}{c} \hat x_{1} \\ \hat x_{2} \end{array} \right]+
\left[ \begin{array}{rrr} 2&1&0\\ 0&1&2 \end{array} \right] \left[
\begin{array}{c}  u_{1} \\ u_{2} \\ u_{3} \end{array} \right].$$
So $\hat W_{AB}=[MB \mid {\hat A}MB]=\left[ \begin{array}{rrrrrr} 2 &1 &0 &-k &0 &k\\
0 &1 &2 &k &0 &-k
\end{array} \right]$ which has rank 2, therefore in this case the \laS\ is also \cC \quad in
accordance with Theorem \ref{redir}.

\subsubsection{Not controllable original---controllable lumped system}
Second, let us consider the control \dE\
$$\frac{d}{dt} \left[ \begin{array}{c}  x_{1} \\ x_{2} \\ x_{3} \end{array} \right]=
\left[ \begin{array}{rrr} -k&k&0 \\ k&-2k&k \\ 0&k&-k \end{array}
\right] \left[ \begin{array}{c}  x_{1} \\ x_{2} \\ x_{3}
\end{array} \right]+ \left[ \begin{array}{rrr} 1&1&0 \\ 1&0&0 \\
1&-1&0
\end{array} \right] \left[ \begin{array}{c}  u_{1} \\ u_{2} \\ u_{3}
\end{array} \right].$$ \label{pelda1} Since in this case the
matrix $W_{AB}= \left[
\begin{array}{rrrrrrrrr} 1 &1 &0 &0 &-k &0 &0 &k^2 &0 \\ 1 &0 &0 &0 &0
&0 &0 &0 &0
\\ 1 &-1 &0 &0 &k &0 &0 &-k^2 &0 \end{array} \right]$
 has only rank 2, thus this system is not \cC . At the same time
the \laS\ is
$$\frac{d}{dt} \left[ \begin{array}{c} \hat x_{1} \\ \hat x_{2} \end{array} \right]=
\left[ \begin{array}{rr} -\frac{k}{2}&\frac{k}{2} \\
\frac{k}{2}&-\frac{k}{2} \end{array} \right] \left[
\begin{array}{c} \hat x_{1} \\ \hat x_{2} \end{array} \right]+
\left[ \begin{array}{rrr} 3 &2 &0 \\ 3 &-2 &0 \end{array} \right]
\left[
\begin{array}{c}  u_{1} \\ u_{2} \\ u_{3} \end{array} \right].$$ Therefore
$\hat W_{AB}=\left[ \begin{array}{rrrrrr} 3 &2 &0 &0 &-2k &0\\ 3
&-2 &0 &0 &2k &0
\end{array} \right]$ and this matrix has rank 2, so the \laS\ is
still \cC\ in this case.

\subsection{Examples for the effect of lumping on observability}
\subsubsection{Observable original---observable lumped system}
First, let us consider the observation system
$$\frac{d}{dt} \left[ \begin{array}{c}  x_{1} \\ x_{2} \\ x_{3} \end{array} \right]=
\left[ \begin{array}{rrr} -k&k&0 \\ k&-2k&k \\ 0&k&-k \end{array} \right]
\left[ \begin{array}{c}  x_{1} \\ x_{2} \\ x_{3} \end{array} \right], \quad
\left[ \begin{array}{c} y_{1} \\ y_{2} \end{array} \right]=
\left[ \begin{array}{rrr} 1&1&0\\ 0&1&1 \end{array} \right]
\left[ \begin{array}{c}  x_{1} \\ x_{2} \\ x_{3} \end{array} \right]$$
i.e. the observation matrix $C=\left[ \begin{array}{rrr} 1&1&0\\ 0&1&1 \end{array} \right].$
So $V_{CA}=[C^{\top} \mid A^{\top}C^{\top} \mid (A^{\top})^2 C^{\top}]=
\left[ \begin{array}{rrrrrr} 1 &0 &0 &k &-k^2 &-2k^2 \\ 1 &1 &-k &-k &3k^2 &3k^2
\\ 0 &1 &k &0 &-2k^2 &-k^2 \end{array} \right]$
which has rank 3, hence this system is \cO\  \quad by Theorem
\ref{megf}. The \laS\ is characterized by $\tilde A= \left[
\begin{array}{rr} -\frac{k}{2}&\frac{k}{2}
\\\frac{k}{2}&-\frac{k}{2} \end{array} \right]$ and $MC^{\top}=\left[ \begin{array}{rr} 3&1 \\
1&3
\end{array} \right]$ by Remark \ref{drem}.
So $[MC^{\top} \mid \tilde AMC^{\top}]=\left[
\begin{array}{rrrr} 3 &1 &-k &k
\\ 1 &3 &k &-k
\end{array} \right]$ which has rank 2, therefore in this case the \laS\ is
\cO\ .
\subsubsection{Not observable original---observable lumped system}
Second, let us consider the observation system:
$$\frac{d}{dt} \left[ \begin{array}{c}  x_{1} \\ x_{2} \\ x_{3} \end{array} \right]=
\left[ \begin{array}{rrr} -k&k&0 \\ k&-2k&k \\ 0&k&-k \end{array}
\right] \left[ \begin{array}{c}  x_{1} \\ x_{2} \\ x_{3}
\end{array} \right], \quad \left[ \begin{array}{c} y_{1} \\ y_{2}
\end{array} \right]= \left[ \begin{array}{rrr} 2&1&0\\ 0&1&2
\end{array} \right] \left[ \begin{array}{c}  x_{1} \\ x_{2} \\
x_{3} \end{array} \right]$$ i.e. now $C=\left[ \begin{array}{rrr}
2&1&0\\ 0&1&2 \end{array} \right]$. \label{pelda2}  Since in this
case the matrix $$V_{CA}= \left[ \begin{array}{cccccc} 2 &0 &-k &k
&k^2 &-k^2 \\ 1 &1 &0 &0 &0 &0
\\ 0 &2 &k &-k &-k^2 &k^2 \end{array} \right]$$
which has rank only 2, thus this system is not \cO\ . At the same
time the \laS\ is $\tilde A=
\left[ \begin{array}{rr} -\frac{k}{2}&\frac{k}{2} \\
\frac{k}{2}&-\frac{k}{2} \end{array} \right]$ and $
MC^{\top}=\left[
\begin{array}{rr} 5&1 \\ 1&5 \end{array} \right].$  Therefore $[MC^{\top} \mid \tilde AMC^{\top}]=
\left[ \begin{array}{rrrr} 5 &1 &-2k &2k \\ 1 &5 &2k &-2k
\end{array} \right]$ and this matrix has rank 2, so the \laS\ is
\cO\ in this case.

\torol{
\section{Appendices}
}%torolvege

%torolvege

%torolvege

%torolvege

\section*{Acknowledgements}
The author acknowledges J. T\'{o}th for giving
important and useful references, and B. Garay and
\'{E}. Gyurkovics for reading
an earlier version of the manuscript.
The research has been partially supported by the
National Scientific Foundation, Hungary, under No. K 63066  and
T047132, and by the French Ministry of Territorial Planning and
Environment (BCRD AP2001).

\section*{Discussion and outlook}

\torol{
 Some of the most important compartmental systems, such a
irreversible catenary, mamillary and circular systems are
symbolically simplified by the method of exact linear lumping.
Transformation of the qualitative properties under lumping are
also traced.

The most important classes
of compartmental systems have been reviewed from the point of view
of symbolic lumpability. Practically interesting lumped systems
mainly arise from numerical calculations, which can be carried out
in all cases without difficulties. We used the sentence "which is
the induced kinetic differential equation of the reaction"
recurrently.

Our aim here is to give explicitly possible lumped compartmental
systems in a few important classes, mainly of symmetric structure
such as: mamillary models, catenary models and circular models.
Some classes can be treated in full generality, some only under
restrictions on the parameters.
}%torolvege

 One of the major questions connected with lumping (with this special
 technique to reduce the number of variables) is: are the qualitative
 properties of the lumped and of the original system connected? This problem
has been investigated in a more general setting in \ref{lump}; here we
add a new statement: suppose we lump a \cC\ system of $n$
compartment. Then the lumped system will also be \cC.
Previously, \cite{F.Gy.c} investigated a similar problem: local
observability and local controllability of reactions.
\cite{E.T.} mainly concentrates on symbolic lumping of compartmental systems.

Possible
further topics are: the effect of nonlinear lumping on local
controllability and observability.

\end{document}